%% file: carbuncle.tex
\documentclass{article}

\usepackage{amsmath}
\usepackage{amsxtra}
\usepackage{amssymb}
\usepackage{amsfonts}
\usepackage{theorem}

\usepackage{graphics}
\usepackage{graphicx}

\usepackage{color}

\newcommand{\mycaption}[1]{\par\refstepcounter{figure}\textsc{Figure} \arabic{figure}. #1}
\newcommand{\defm}[1]{\textit{#1}}

\title{Carbuncles as self-similar entropy solutions}
\author{Volker Elling}
\date{}

\begin{document}

\maketitle

\parindent=0cm
\parskip=\baselineskip

\begin{abstract}
\noindent Numerical approximations of shock waves sometimes suffer from instabilities called \defm{carbuncles}. Techniques for suppressing carbuncles are trial-and-error and lack in reliability and generality, partly because theoretical knowledge about carbuncles is equally unsatisfactory. It is not known which numerical schemes are affected in which circumstances, what causes carbuncles to appear and whether carbuncles are purely numerical artifacts or rather features of a continuum equation or model. This work presents evidence towards the latter: it is conjectured that carbuncles are a special class of non-physical entropy solutions. Using a new technique for triggering a single carbuncle, its structure is computed in detail in similarity coordinates.
\end{abstract}

\section{Introduction}

\begin{figure}
\parbox[t]{.45\textwidth}{%
\includegraphics[width=\linewidth]{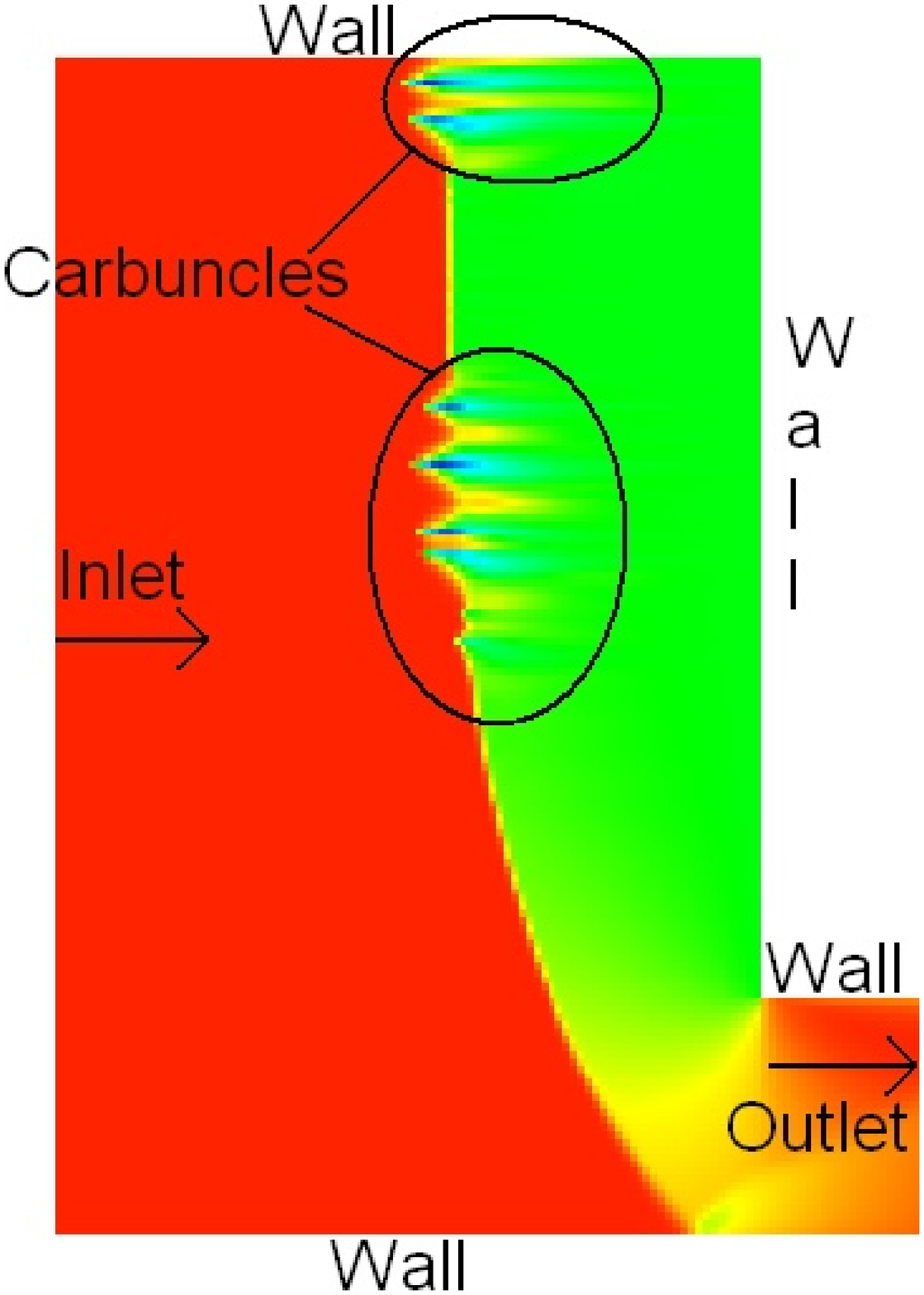}%
\mycaption{Curved shock with carbuncles in front of a forward-facing step}
\label{fig:forwardstep}
}
\hspace{1cm}%
\parbox[t]{.5\textwidth}{%
\includegraphics[width=\linewidth]{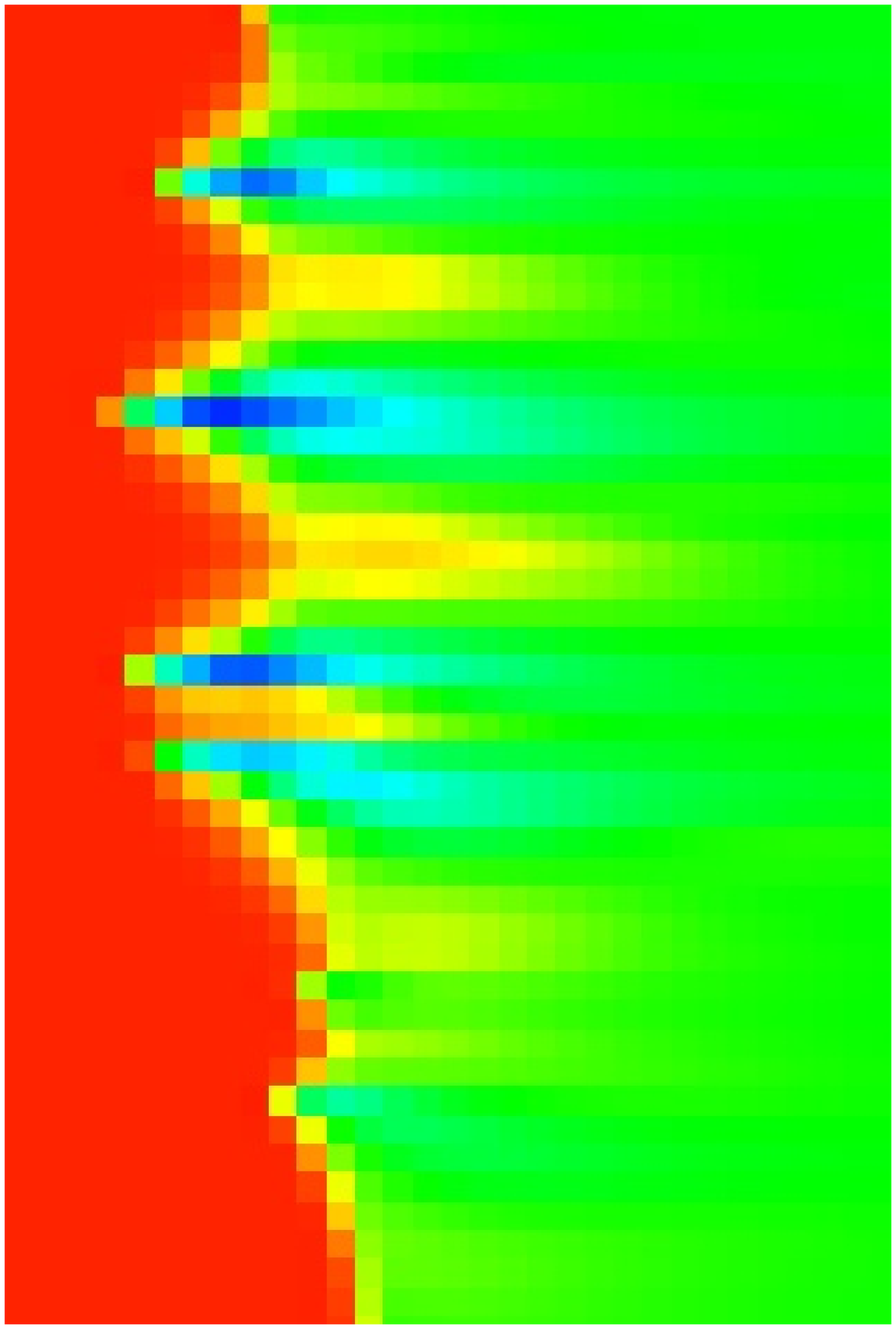}%
\mycaption{Closeup of carbuncles in Figure \ref{fig:forwardstep}}
\label{fig:forwardstepdetail}
}
\vspace{5mm}

\includegraphics[width=\textwidth]{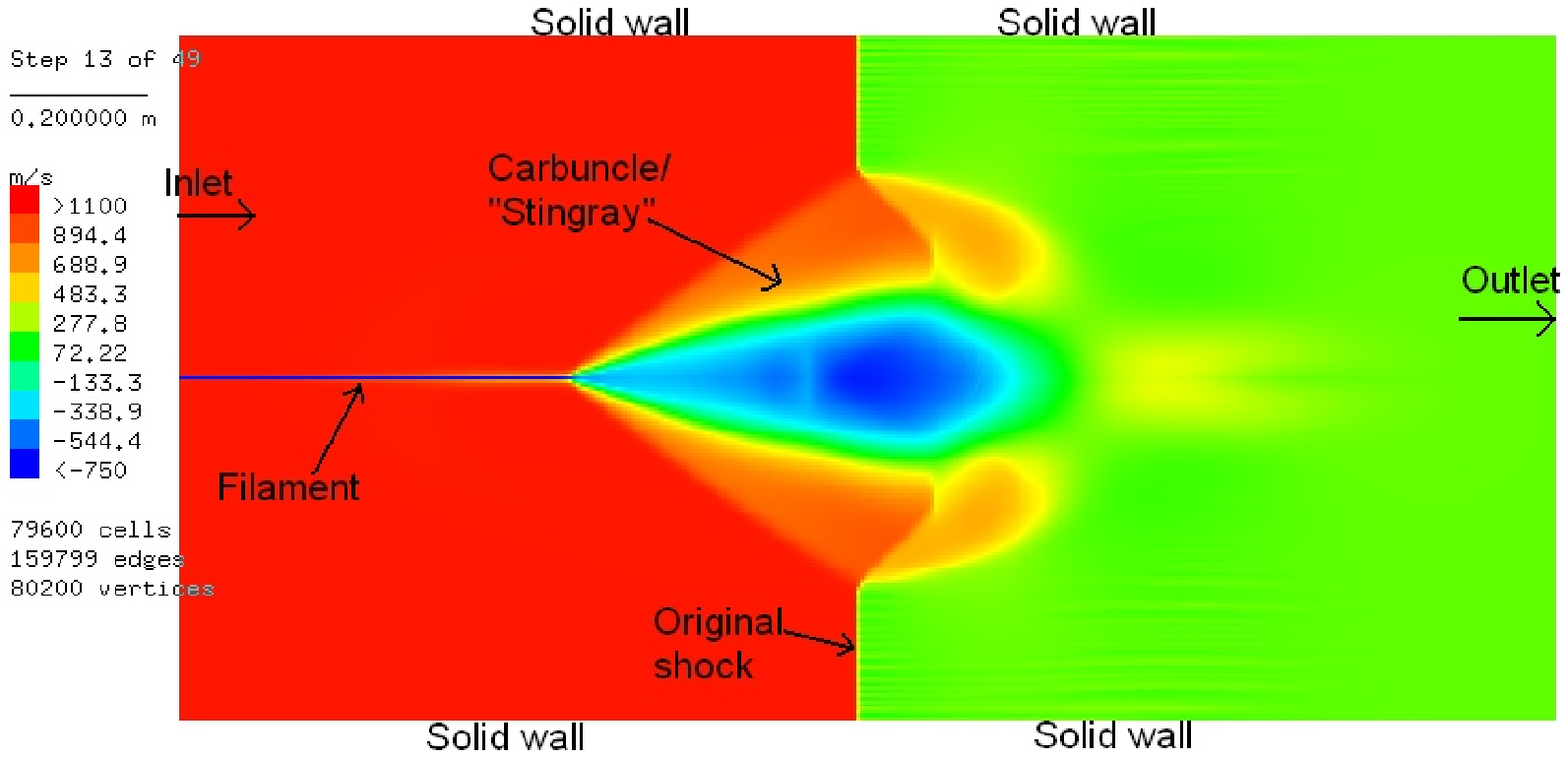}%
\mycaption{A stingray-like structure resembling carbuncles}
\label{fig:steady-carbuncle}
\end{figure}

Numerical calculations of shock waves in inviscid compressible flow sometimes suffer from instabilities, which were christened \defm{carbuncles}
in \cite{peery-imlay}.
Figure \ref{fig:forwardstep} presents an example: a shock wave in front of a forward-facing step, normally expected to be smooth, develops several small structures called ``carbuncles'' in the literature. Figure \ref{fig:forwardstepdetail} shows these carbuncles in more detail. For other examples of carbuncle phenomena see \cite{peery-imlay}, \cite{dumbser-moschetta-gressier} or \cite[Figure 5]{quirk}. 

There is no clear rule for predicting which numerical methods are likely to be affected or in which circumstances carbuncles might appear (a summary of various proposals and explanations can be found in \cite{dumbser-moschetta-gressier}). Figure \ref{fig:forwardstep} displays two characteristic situations: the upper group of carbuncles is located near a boundary where numerical boundary layers act; the lower group appears in a region where the shock is almost but not exactly parallel to the grid edges. Some authors single out the Roe scheme \cite{roe} as particularly prone to carbuncles whereas others report bad experiences with the Godunov method \cite{godunov}. It appears that the theoretical consensus on carbuncles is as volatile as the carbuncles themselves. However, recently researcher have focused on instability of discrete plane shock waves as an ingredient for the genesis of carbuncles (see e.g.\ \cite{dumbser-moschetta-gressier,robinet-gressier-casalis-moschetta,aiso}).

As a consequence of the theoretical confusion, there are no safe recipes for eliminating carbuncles. Although numerical viscosity in various forms can suppress the carbuncle phenomenon, the exact amount required is not known. Since numerical viscosity smears discontinuities and decreases overall accuracy, it is desirable to have precise theoretical criteria for carbuncle avoidance.

Due to the nature of carbuncles as phenomena without precise mathematical definition, it is not clear whether there is unanimous agreement in the scientific community on which phenomena can be classified as carbuncles. Numerical calculations of discontinuities in gas flow can suffer from many other problems, such as linear instability due to excessive time steps or poor choices for operator discretization, nonlinear phenomena in the presence of large total variation, Kelvin-Helmholtz instabilities, expansion shocks etc. We hope the reader will agree that the quasi-continuum carbuncle patterns computed in standard (Figure \ref{fig:steady-carbuncle}) and similarity (Figure \ref{fig:carbuncle-vx-detail}) coordinates in this paper correspond to to the small discrete carbuncles shown in Figures \ref{fig:forwardstep} and \ref{fig:forwardstepdetail} or in the references above.

All results displayed in this paper were done for the 2D isentropic Euler equations with $\gamma:=7/5$; however, analogous results can be generated for other $\gamma$ and for the nonisentropic Euler equations as well. In all diagrams, flow into the domains is Mach 3.0 horizontally from left to right. Unless otherwise noted, all pictures show the horizontal component of velocity.

\section{Carbuncles in standard coordinates}

In this section we describe how carbuncles can be triggered in a large variety of numerical schemes in standard coordinates (here, ''standard'' means coordinates $(t,x)$ as compared to similarity coordinates $(t,\frac{x}{t})$.) We start with a steady vertical plane shock wave computation in a uniform Cartesian grid. The upper and lower domain boundaries are solid walls with slip boundary condition; fluxes across the left resp.\ right domain boundary are computed with pre- resp.\ post-shock state in ``pseudo-cells'' on the outside of the domain boundaries. In this setting, the exact plane shock wave is reproduced faithfully by discrete shocks.

Now the initial and boundary data is modified as follows: we choose an edge in the center of the left boundary; in its pseudo-cell and in all cells in a horizontal one-cell-high ``filament'' from this edge to the shock, the horizontal velocity is set to zero. This filament immediately triggers a structure (see Figure \ref{fig:steady-carbuncle}) that grows steadily (and eventually reaches the domain boundaries). This structure, somewhat reminiscent of a stingray, is very similar to the small discrete carbuncles shown in Figure \ref{fig:forwardstepdetail}.

Figure \ref{fig:steady-carbuncle} has been computed with the Godunov scheme \cite{godunov}; similar carbuncles can be observed in the local Lax-Friedrichs \cite{shu-osher-llf} or the Osher-Solomon \cite{solomon-osher}  schemes. Higher-order schemes are likewise affected (this is not particularly surprising because the creation of the carbuncle takes place in nonsmooth or even discontinuous regions where higher-order schemes revert to first order). Moreover, adding generous amounts of numerical viscosity anywhere except near the filament and the carbuncle tip does not destroy the carbuncle. This suggests that the tip is the unphysical part of the carbuncle (unless of course carbuncles are physically meaningful).

The reader may suspect that a pronounced disturbance like the carbuncle in Figure \ref{fig:steady-carbuncle} is inevitable given the reduced mass and momentum transport in the filament. However, the filament is too thin to support this reasoning. When the grid is refined and the filament height (and hence the $L^1$ norm perturbation to the initial and boundary data) decreases, the carbuncle does not diminuish; it always grows until it pushes the confines of the domain.

\section{Carbuncles in similarity coordinates}

\begin{figure}
\parbox[t]{.43\textwidth}{%
\includegraphics[width=\linewidth]{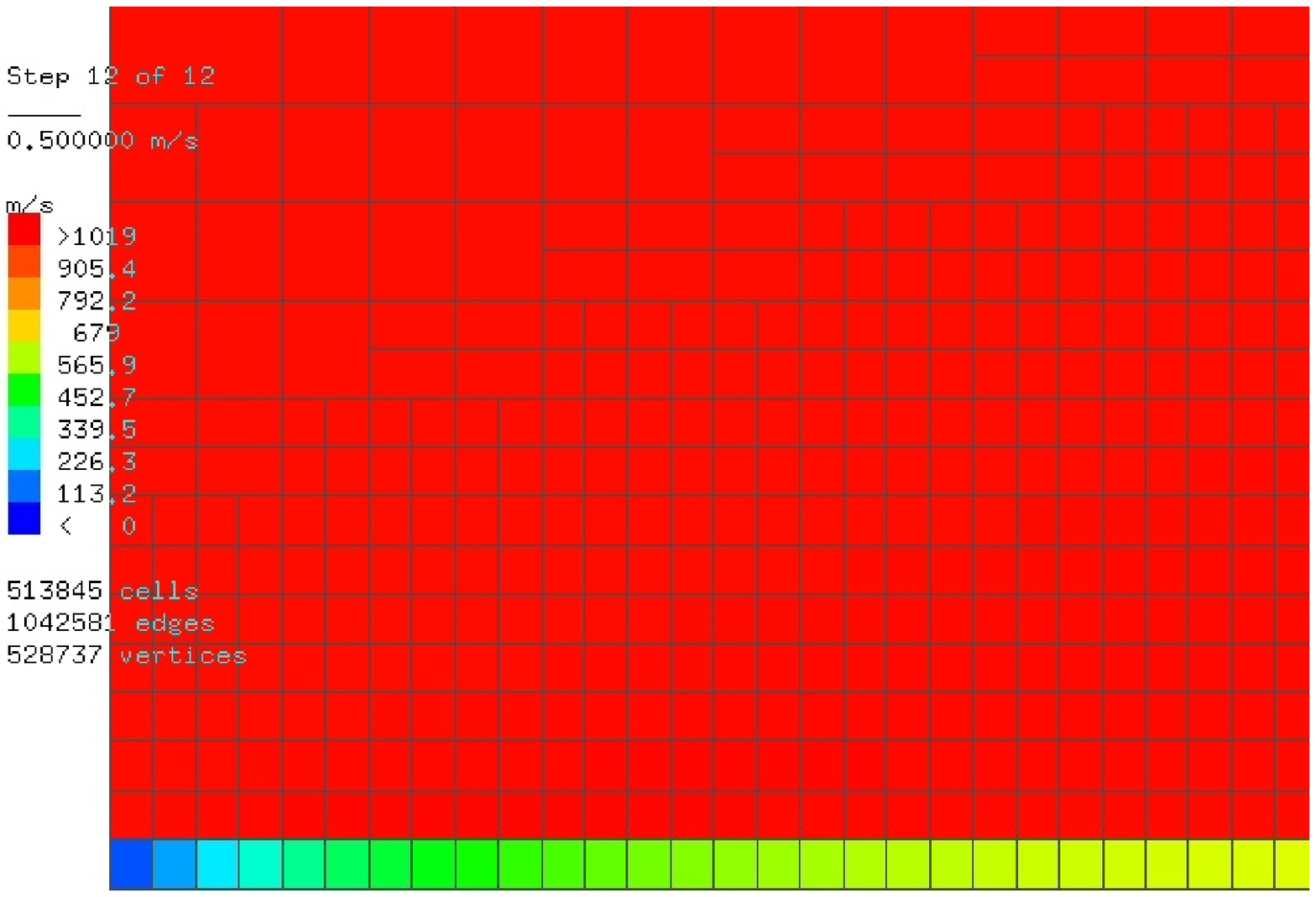}%
\mycaption{Filament emanating from the lowest edge in the left boundary}
\label{fig:trigger}
}
\parbox[t]{.55\textwidth}{%
\includegraphics[width=\linewidth]{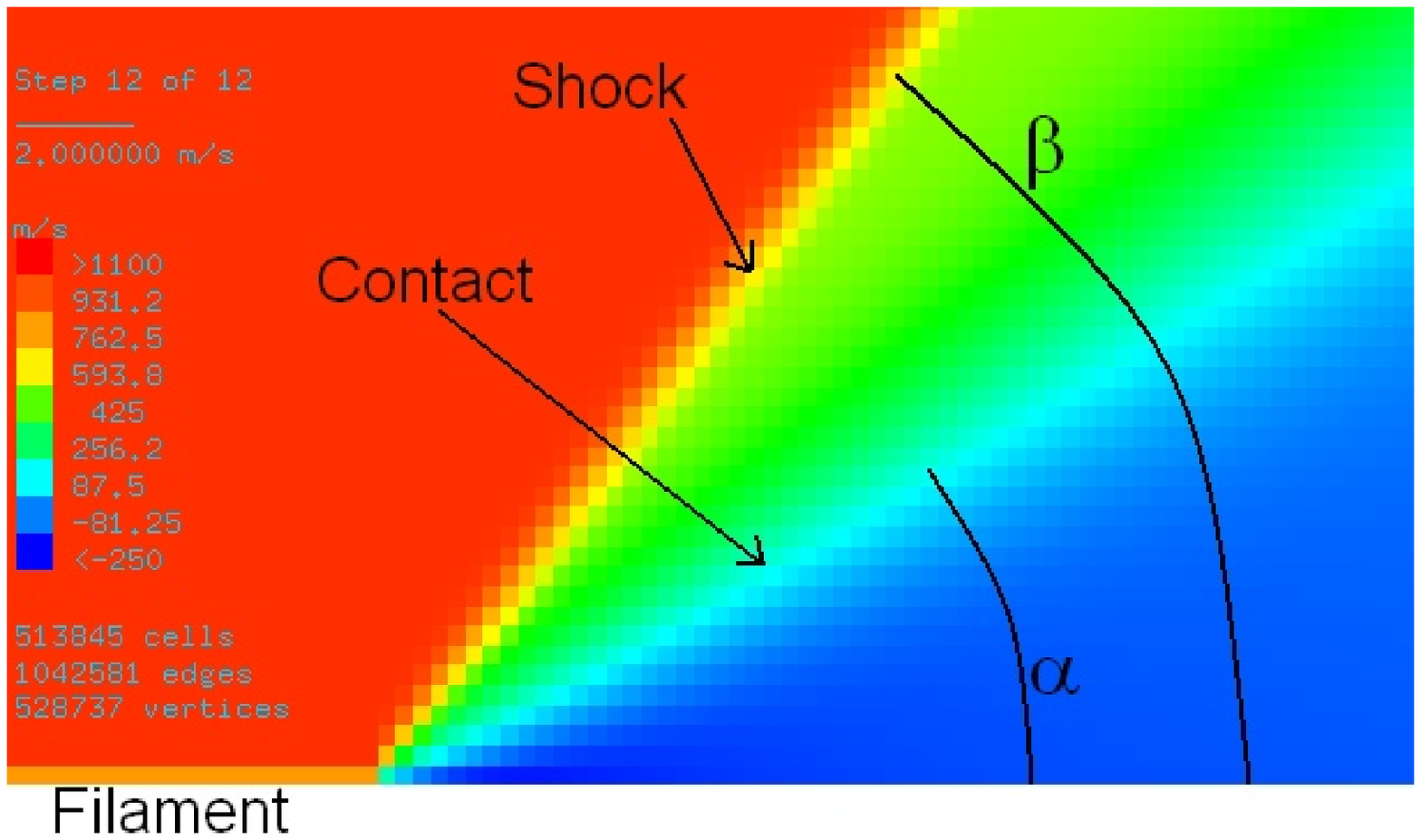}%
\mycaption{Filament and carbuncle tip (contact smeared due to poor grid alignment)}
\label{fig:tip}
}
\vspace{5mm}

\includegraphics[width=\textwidth]{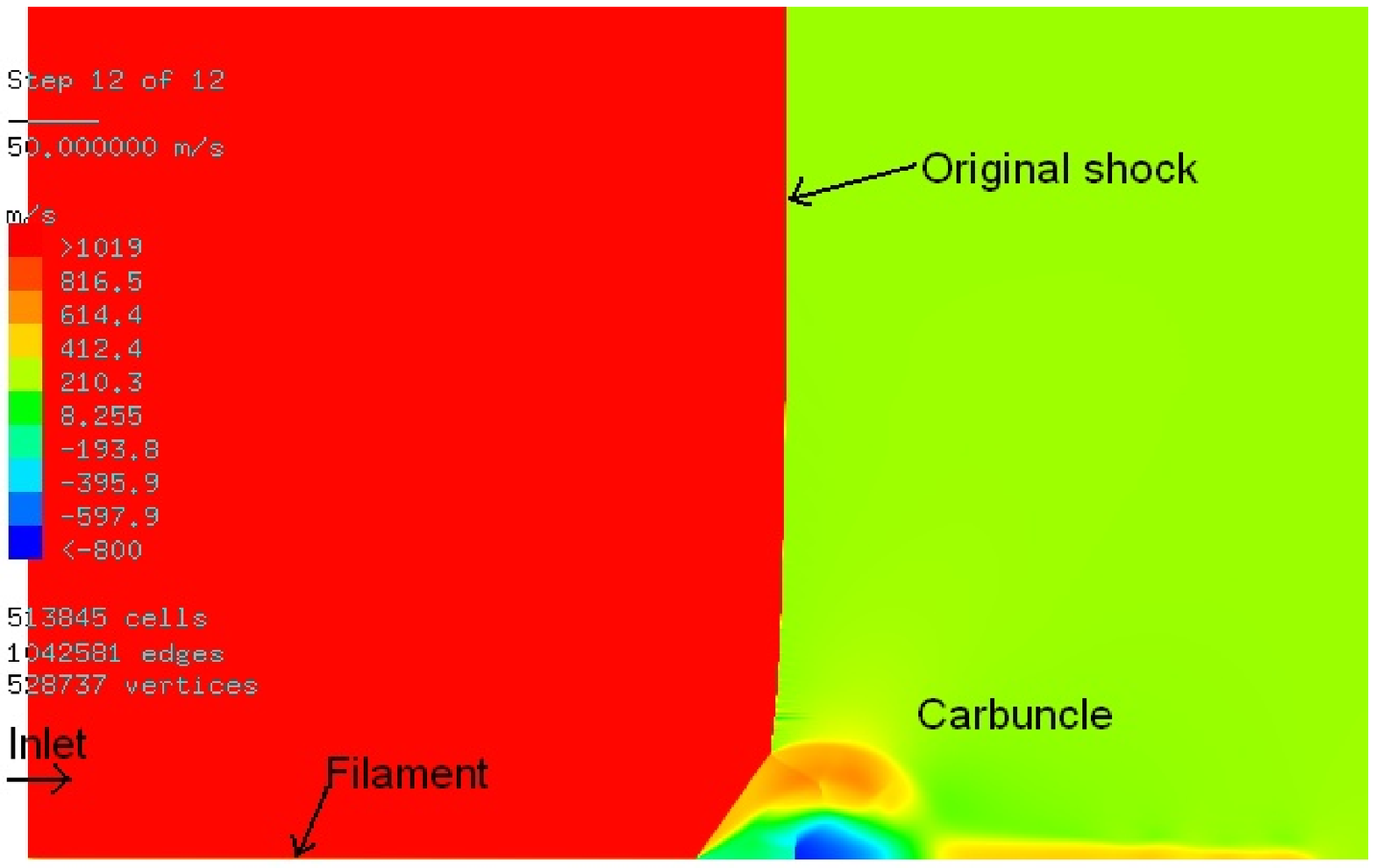}%
\mycaption{Carbuncle in a plane shock (in similarity coordinates)}
\label{fig:carbuncle-vx-whole}
\end{figure}

\begin{figure}
\includegraphics[width=\textwidth]{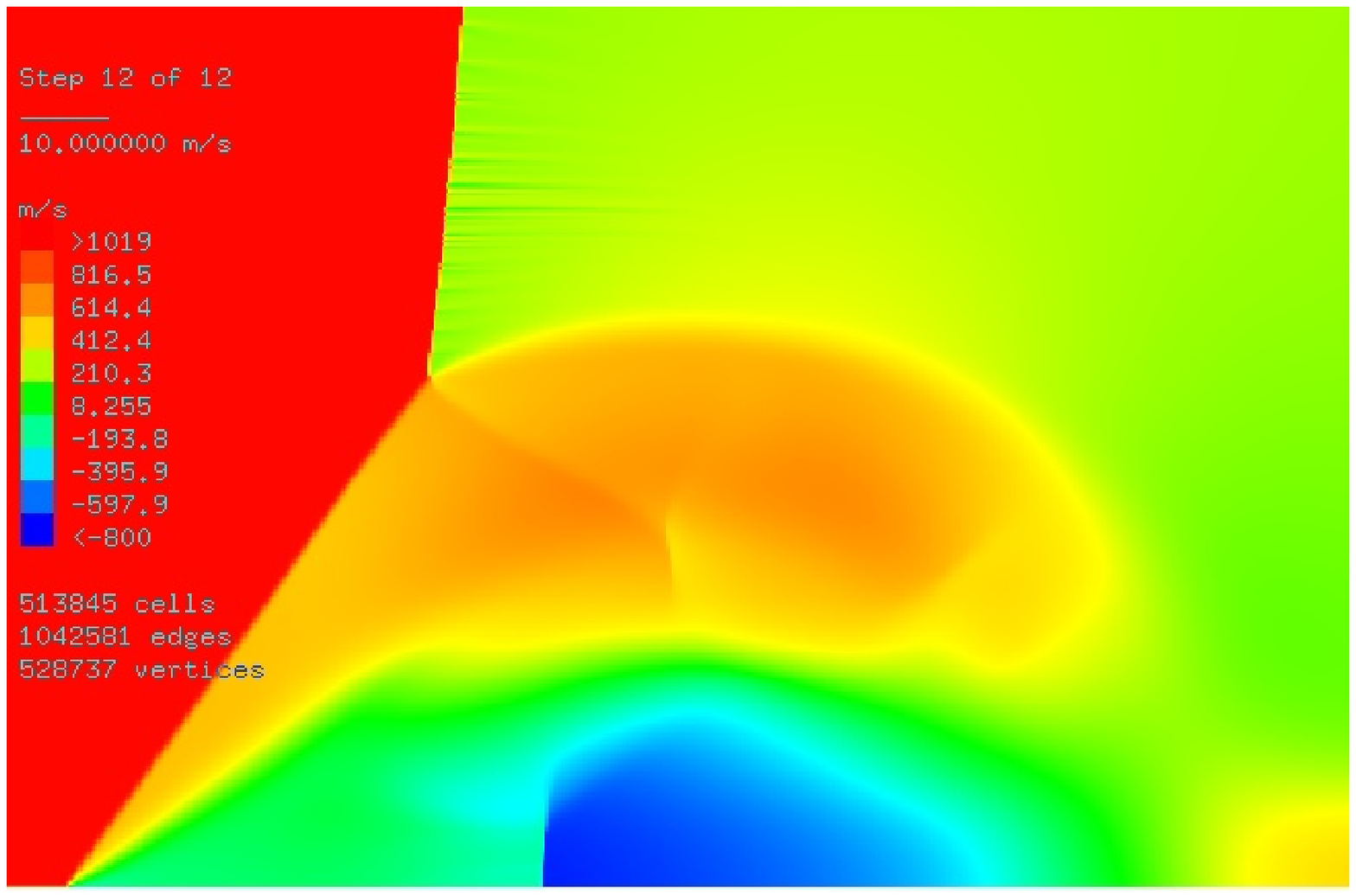}%
\mycaption{Detailed self-similar structure of carbuncle in Figure \ref{fig:carbuncle-vx-whole}}
\label{fig:carbuncle-vx-detail}
\vspace{5mm}

\includegraphics[width=\textwidth]{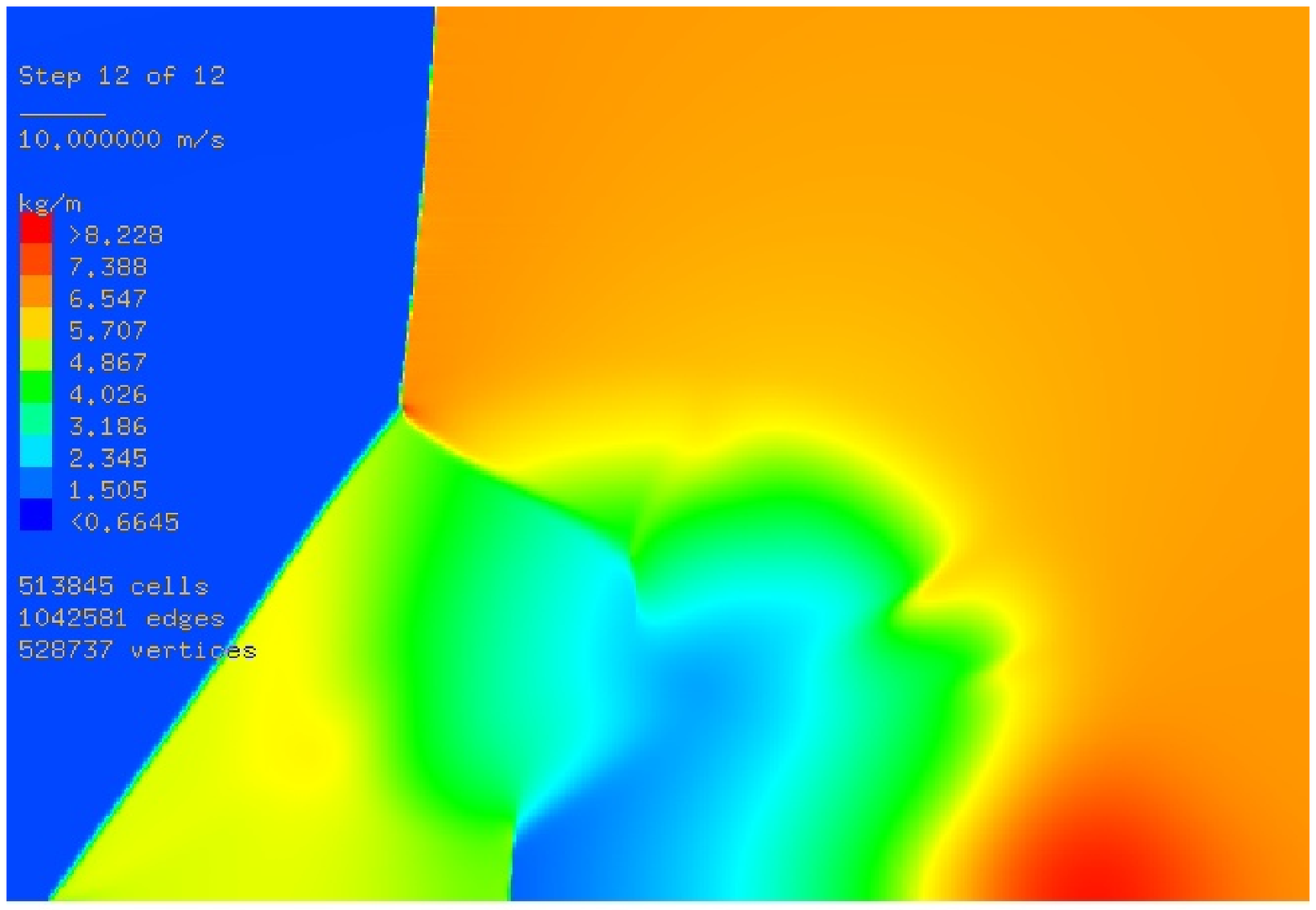}%
\mycaption{Density field for Figure \ref{fig:carbuncle-vx-detail}}
\label{fig:carbuncle-rho-detail}
\end{figure}

The new method for triggering carbuncles is most effective in standard coordinates where it can be successfully applied to many numerical methods. It is observed that the carbuncle pattern grows at a constant velocity and that its inner structure settles in after some time. This suggests that carbuncles are self-similar flow patterns. Moreover, the steadily growing carbuncles eventually reach the boundary of any spatial domain. Therefore, for computing their structure in detail similarity coordinates $(t,\frac{\vec x}{t})$ are more appropriate.

A major drawback is that similarity coordinates add dissipation because grid cells move along the shock and smoothen the solution by averaging. For most numerical schemes, this decreases carbuncles to a point where it cannot be said for sure whether carbuncle size is independent of filament size. More importantly, it is not possible to create carbuncles of sufficient size to compute significant detail. However, the Godunov scheme does produce robust and large carbuncles in our setting.

Here we compute on a Cartesian grid with adaptive refinement to achieve better resolution. Similarity coordinates are treated by the moving-edge modifications discussed in \cite{elling-diplom,phd-thesis}. To save time, only the upper halfplane of the obviously symmetric carbuncles is computed. The trigger mechanism remains the same: in a one-cell-high filament at the lower domain boundary, horizontal pre-shock velocity is reduced to $0$. In similarity coordinates the resulting contact discontinuity is not preserved exactly; additional waves propagate upwards and the filament weakens noticeably (see Figure \ref{fig:trigger}). However, this unelegant sideeffect is not too large for our purposes. With the filament located at the lower boundary, it has a lot of similarity with numerical boundary layers. Indeed, in the experience of the author these are the most frequent cause of carbuncle-like phenomena.

In this new setting, carbuncle patterns can be computed in great detail (see Figure \ref{fig:carbuncle-vx-whole} for overview and Figures \ref{fig:carbuncle-vx-detail} and \ref{fig:carbuncle-rho-detail} for more detail). Better resolution could be achieved by using higher-order modifications, but here we choose the Godunov scheme because of its theoretical appeal (it manifestly satisfies all discrete entropy inequalities).

\section{Relation to entropy solution non-uniqueness}

\begin{figure}
\input{example.pstex_t}
\mycaption{Probably unphysical steady self-similar entropy solution}
\label{fig:nuq-initial}
\vspace{5mm}

\includegraphics[height=.4\textheight]{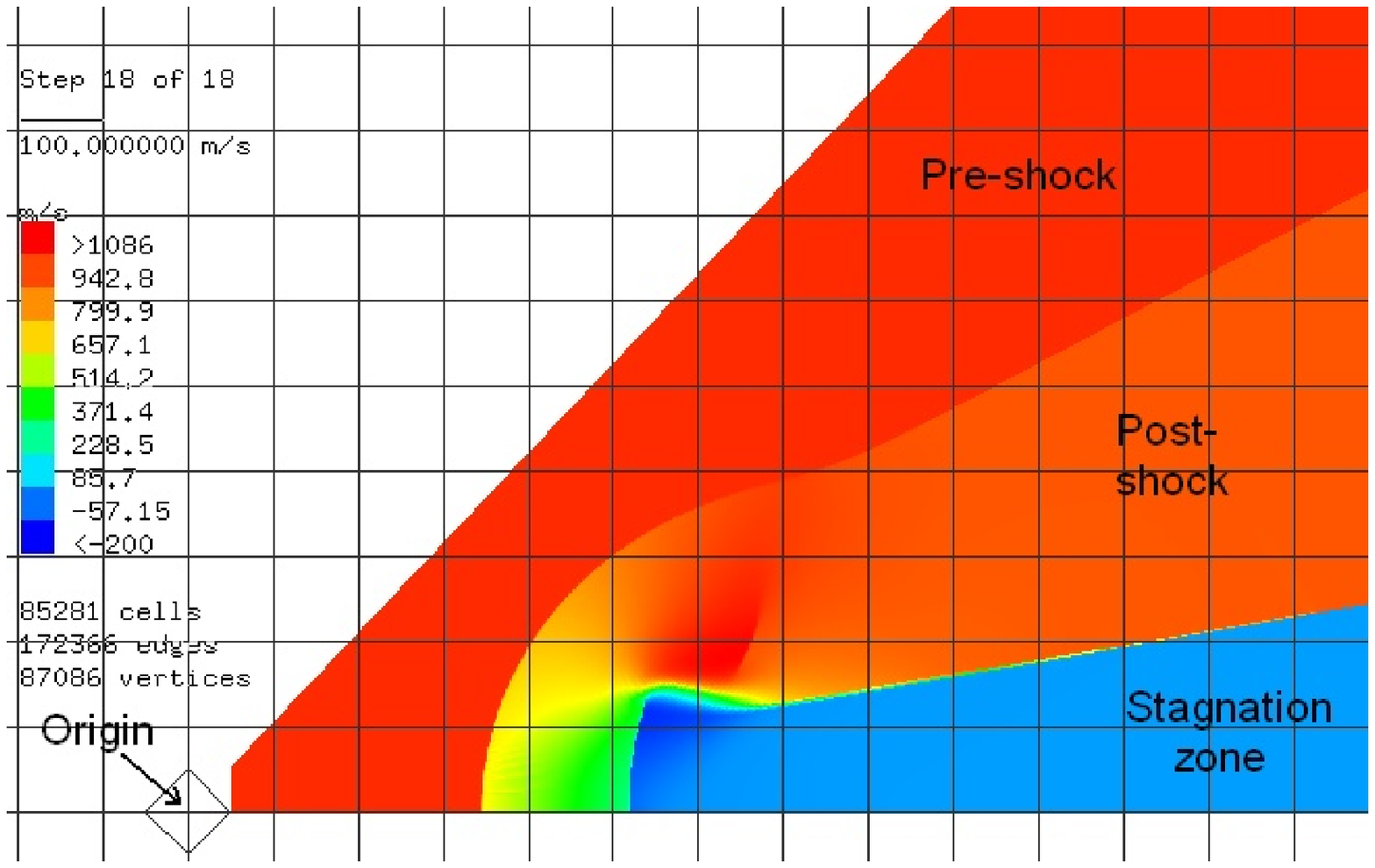}%
\mycaption{A second self-similar (but unsteady) entropy solution, for Figure \ref{fig:nuq-initial} as initial data}
\label{fig:nuq-physical}
\end{figure}

In \cite{elling-nuq-journal} a numerical example is presented that suggests that entropy solutions of the Cauchy problem for the compressible Euler equations are not always unique. More precisely, Figure \ref{fig:nuq-initial} constitutes an exact steady self-similar entropy solution; however, using it as initial data, numerical schemes produce another solution, Figure \ref{fig:nuq-physical}, which appears to be self-similar and (being computed by the Godunov scheme) an entropy solution, but clearly unsteady.

The initial data in this example was motivated by research on the problem of supersonic flow onto a solid wedge. For small wedge opening angles ($\alpha$ in Figure \ref{fig:nuq-initial}), there are two possible solutions, a strong and a (comparatively) weak shock (which was chosen for computing Figure \ref{fig:nuq-physical}). Replacing the solid wedge by gas with post-shock density and zero velocity yields Figure \ref{fig:nuq-initial}.

In Figure \ref{fig:tip} the tip of the carbuncle in Figure \ref{fig:carbuncle-vx-whole} is shown in detail. Clearly it resembles Figure \ref{fig:nuq-initial} (the post-shock area in Figure \ref{fig:tip} is supersonic which means the shock corresponds to the \emph{weak} shock in the wedge problem). This analogy suggests that carbuncles are another manifestation of the same problem: two different entropy solutions for identical data. In this case, the \emph{theoretical} solution (the plane shock wave) is known to be physically correct as well as a vanishing viscosity and hence entropy solution; the carbuncle appears to be a second entropy solution which is probably unphysical without the presence of the filament.

This paper is not the first to conjecture that carbuncles are phenomena of the continuum equations rather than purely numerical artefacts. Other references are \cite{robinet-gressier-casalis-moschetta} where a linear stability analysis of plane shocks in the exact Euler equations seems to reveal a new unstable mode related to carbuncles (\cite{coulombel-benzoni-gavage-serre} have raised objections to this analysis), as well as \cite{morton-roe} and \cite{dumbser-moschetta-gressier} both of which explicitly state suspicions that carbuncles and other anomalies are entropy solutions.

It should be emphasized that, although the carbuncle in Figure \ref{fig:carbuncle-vx-whole} does not appear to vanish as the grid is refined, it shrinks and grows and does not seem to converge to a single pattern either. Indeed, if the comparison to the nonuniqueness example is appropriate, then there may be a parametrized family of carbuncle patterns rather than a single one. There are at least two parameters that are currently not under our control: the tip angle $\beta$ and the location of the carbuncle tip. These are determined in some unknown and probably rather sensitive way by the numerical method, the filament and the inner structure of the carbuncle pattern. If one or both parameters are free (i.e.\ if there is a different carbuncle pattern for each choice of parameter in some interval), the carbuncle calculation should not be expected to converge in general.

More generally, we should stress that even on a single fixed grid Figure \ref{fig:carbuncle-vx-whole} is rather unstable; small modifications to the numerical scheme can have drastic effects. For example, the Godunov scheme requires solving a nonlinear equation for every edge in every time step; solving this equation with less than machine accuracy can either trigger oscillations that abort the computation or produce additional smaller carbuncles along the shock that grow to spoil the overall pattern.

It is possible to stabilize the calculation to some degree by adding numerical viscosity everywhere except near the lower boundary; this suppresses the appearance of additional carbuncles, but smears the carbuncle pattern even more.

\section{Conclusions and open problems}

This article presents a new technique for triggering carbuncles in shock waves. Its main advantage is reliability; previous carbuncle computations were rather precarious and sensitive to changes in setting and numerical methods. There are two immediate benefits: first, it becomes clear that a rather large variety of numerical methods is affected, more than may have been suspected previously. Second, carbuncle structure can be computed in fine detail. This bears fruit immediately: as discussed in the previous section, Figure \ref{fig:tip} suggests a connection between the carbuncle phenomenon and the hypothetical entropy solution nonuniqueness example (see above). Moreover, it can be verified that carbuncles are self-similar structures. Numerical experiments in which carbuncles appear to be steady rather than self-similar most likely exhibit equilibrum states in which further carbuncle growth would weaken numerical effects that stabilize the carbuncle tip.

The increased understanding of the nature of carbuncles should facilitate the development of better criteria for predicting the appearance of carbuncles and of numerical methods that avoid carbuncles at a minimal loss of accuracy. Most importantly, there is clear evidence that carbuncles are entropy solutions of the compressible Euler equations rather than non-entropy solutions or numerical artifacts. If carbuncle phenomena are indeed caused by a narrow family of continuum functions, suppressing them might be easier than previously expected.

Nevertheless, this work opens some new questions and leaves several unsolved problems. The most central ones are closely related: first, it is not known how many parameters index the family of carbuncles. Drawing on the analogy to the nonuniqueness example, there could be two parameters, tip angle $\beta$ and tip location. It can be observed that these parameters change as the grid is refined which suggests that at least one of these parameters is free; there is no indication whether the other one is free or rigid.

Moreover, it is desirable to have techniques for \emph{controlling} the free parameters during computations. Using this control, the numerical scheme could be forced to \emph{converge} to a single carbuncle. This would add substantial evidence for the hypothesis that carbuncles are entropy solutions: according to the Lax-Wendroff theorem (see \cite{lax-wendroff,kroener-rokyta-wierse,godlewski-raviart,elling-lax-wendroff-journal}), a convergent sequence of approximations generated by the Godunov scheme must have an entropy solution as limit. Moreover, with tight control of the free parameters, carbuncle structure could be computed faster, more precisely and more reliably. 

The internal structure in Figure \ref{fig:carbuncle-vx-detail} is quite smeared, so it will change significantly as the grid becomes finer (even in cases where the free parameters happen/are controlled to remain unchanged). Besides lack of control over free parameters, another obstacle towards computing the structure accurately is the strong smearing of contact discontinuities that are not grid-aligned.

Even in the poor resolution of Figure \ref{fig:carbuncle-vx-detail} it is obvious that the carbuncle structure is rather complex. This might in part be due to the rather large pre-shock Mach number of 3.0 which was chosen for the pronounced carbuncle it causes. The carbuncle structure may be simpler for Mach numbers closer to 1, but for such values the author was not able to trigger a carbuncle in similarity coordinates. In standard coordinates, noticeable carbuncles start to appear at about Mach 1.4. It is not clear whether this number is somehow significant or merely a limitation of the filament trigger mechanism. The carbuncles grow very slowly near Mach 1.4, so they may simply escape observation for smaller Mach numbers.

In Figures \ref{fig:nuq-initial} and \ref{fig:nuq-physical} there appear to be two different entropy solutions. This by itself does not rule out that there are only a few isolated entropy solutions; in this case the entropy condition would still be rather useful. However, if the tip shock opening angle is a free parameter that can be shifted up to $180^\circ$ (this requires a strong shock at the carbuncle tip and causes the angle $\alpha$ between the contacts to decrease to $0^\circ$), there may be a continuous (in $L^1$) family of carbuncle patterns, containing the physically correct plane shock wave as a limit case --- such a continuum of entropy solutions for the same data would seriously reduce the value of the entropy condition. Numerical calculations can be used to verify these ideas once we find techniques to control the tip angle $\beta$ (this is another motivation for trying to control the free parameters). However, due to this theoretical relevance of carbuncles, it would be desirable to prove existence of such patterns mathematically rather than numerically. Due to the complex structure of the carbuncle patterns and the intricacies of the Euler equations, this appears to be extremely difficult.

\bibliographystyle{amsalpha}
\bibliography{../../pmeyer/elling}

\end{document}

%% file: example.pstex_t
\begin{picture}(0,0)%
\includegraphics{example.pstex}%
\end{picture}%
\setlength{\unitlength}{3947sp}%
\begingroup\makeatletter\ifx\SetFigFont\undefined%
\gdef\SetFigFont#1#2#3#4#5{%
  \reset@font\fontsize{#1}{#2pt}%
  \fontfamily{#3}\fontseries{#4}\fontshape{#5}%
  \selectfont}%
\fi\endgroup%
\begin{picture}(5371,3675)(214,-4606)
\put(226,-2011){\makebox(0,0)[lb]{\smash{{\SetFigFont{10}{12.0}{\rmdefault}{\mddefault}{\updefault}{\color[rgb]{0,0,0}Inflow (supersonic)}%
}}}}
\put(2701,-1111){\makebox(0,0)[lb]{\smash{{\SetFigFont{10}{12.0}{\familydefault}{\mddefault}{\updefault}{\color[rgb]{0,0,0}Numerical}%
}}}}
\put(2701,-1261){\makebox(0,0)[lb]{\smash{{\SetFigFont{10}{12.0}{\familydefault}{\mddefault}{\updefault}{\color[rgb]{0,0,0}domain}%
}}}}
\put(4201,-2536){\makebox(0,0)[lb]{\smash{{\SetFigFont{10}{12.0}{\rmdefault}{\mddefault}{\updefault}{\color[rgb]{0,0,0}Stagnation zone}%
}}}}
\put(4201,-1036){\makebox(0,0)[lb]{\smash{{\SetFigFont{10}{12.0}{\rmdefault}{\mddefault}{\updefault}{\color[rgb]{0,0,0}Weak shock}%
}}}}
\put(1801,-4561){\makebox(0,0)[lb]{\smash{{\SetFigFont{10}{12.0}{\rmdefault}{\mddefault}{\updefault}{\color[rgb]{0,0,0}$y$}%
}}}}
\put(976,-2686){\makebox(0,0)[lb]{\smash{{\SetFigFont{10}{12.0}{\rmdefault}{\mddefault}{\updefault}{\color[rgb]{0,0,0}$x$}%
}}}}
\put(3151,-2086){\makebox(0,0)[lb]{\smash{{\SetFigFont{12}{14.4}{\rmdefault}{\mddefault}{\updefault}{\color[rgb]{0,0,0}$\beta$}%
}}}}
\put(4201,-2236){\makebox(0,0)[lb]{\smash{{\SetFigFont{10}{12.0}{\rmdefault}{\mddefault}{\updefault}{\color[rgb]{0,0,0}Contact discontinuity}%
}}}}
\put(2926,-2686){\makebox(0,0)[lb]{\smash{{\SetFigFont{10}{12.0}{\rmdefault}{\mddefault}{\updefault}{\color[rgb]{0,0,0}$\alpha$}%
}}}}
\end{picture}%